\documentclass[a4paper, BCOR=0.0mm, DIV=calc]{scrartcl}
\usepackage[T1]{fontenc}
\usepackage[utf8]{inputenc}
\usepackage[english]{babel}
\usepackage[osf,sc]{mathpazo}
\usepackage{textcomp}
\usepackage{microtype}
\usepackage{graphicx, color}
\usepackage{amsmath, amssymb, amsfonts}
\usepackage{booktabs}
\usepackage{nicefrac}
\usepackage{wasysym,pifont}
\usepackage{tikz}
\usepackage{graphicx}
\usepackage{float}
\usepackage{wrapfig}
\usetikzlibrary{calc,intersections}
\KOMAoptions{twoside=false, twocolumn=false, headinclude=false, footinclude=false, mpinclude=false, pagesize=auto}
\recalctypearea

\setkomafont{section}{\mdseries\scshape\Large}
\addtokomafont{subsection}{\mdseries\scshape}
\setkomafont{title}{\mdseries\scshape}

\begin{document}
\title{Euler and the Multiplication Formula for the $\Gamma$-Function.}
\author{Alexander Aycock}
\date{ }
\maketitle

\begin{abstract}
We show that an apparently overlooked result of Euler from \cite{E421} is essentially 
equivalent to the general multiplication formula for the $\Gamma$-function that was 
proven by Gauss in \cite{Ga28}.
\end{abstract}

\section{Introduction}
The interpolation of the factorial by the $\Gamma$-function was found nearly simultaneously by Bernoulli \cite{Be29} and Euler in $1729$ \cite{E19} and is without any doubt one of the most important functions in mathematics. Most of its basic properties were discovered by Euler, who also gave the definition that is nowadays often used to introduce the function in $\S 7$ of \cite{E675}
\begin{equation*}
\Gamma(x) := \int\limits_{0}^{\infty} e^{-t}t^{x-1}dt \quad \text{for} \quad \text{Re}(x)>0.
\end{equation*}
However, a good understanding of $\Gamma$ as a meromorphic function of 
its argument could of course only be achieved after Gauss in the 19th century; the now universally adopted notation
stems from Legrendre \cite{Le09}.\\
One of the fundamental properties of the $\Gamma$-function is so-called {\em multiplication formula} that 
reads, in the modern notation
\begin{equation}
\Gamma \left(\dfrac{x}{n}\right)\Gamma \left(\dfrac{x+1}{n}\right)\cdots \Gamma \left(\dfrac{x+n-1}{n}\right)= \dfrac{(2\pi)^{\frac{n-1}{2}}}{n^{x-\frac{1}{2}}}\cdot \Gamma(x).
\end{equation}
For $n=2$ one obtains the {\em duplication formula} that is usually ascribed to Legendre \cite{Le26}.\\[2mm]
The multiplication formula was first proven rigourously by Gauss in  his influential paper \cite{Ga28} 
on the hypergeometric series, in which he also gives a complete account of the factorial function $\Pi(x):=\Gamma(x+1)=x!$.  
Gauss cited Euler's results very often, but apparently he was not aware of the lesser-known paper \cite{E421} of Euler.
In that paper Euler presents a formula that is essentially equivalent (1), as we will explain now.

\subsection{The function $\left(\frac{p}{q}\right)$}

In $\S 3$ of \cite{E321} and $\S 44$ of \cite{E421}, Euler studies properties of the function 
\begin{equation*}
\left(\frac{p}{q}\right):= \int\limits_{0}^{1} \dfrac{x^{p-1}dx}{(1-x^n)^{\frac{n-q}{n}}}.
\end{equation*}
In his notation the variable $n$ is left implicit, and Euler shows the nice symmetry property
\[ \left(\frac{p}{q}\right)=\left(\frac{q}{p}\right) .\]
Of course, by the substitution $x^n=y$ this function is just the Beta-function in disguise: 
\begin{equation}
\left(\frac{p}{q}\right) = \frac{1}{n} \int\limits_{0}^1 y^{\frac{p}{n}-1}dy(1-y)^{\frac{q}{n}-1}=\frac{1}{n} B \left(\frac{p}{n},\frac{q}{n}\right),
\end{equation}
where the Beta-function is defined as
\begin{equation*}
B(x,y)= \int\limits_{0}^1 t^{x-1}dt(1-t)^{y-1} \quad \text{for} \quad \text{Re}(x),\text{Re}(y) > 0.
\end{equation*}
Euler implicitly assumes $p$ and $q$ to be natural numbers,  but
this restriction is of course not necessary.\\
One of the early discoveries of Euler \cite{E19} was that the Beta-integral reduces to a product of $\Gamma$-factors:
\begin{equation*}
B(x,y) = \dfrac{\Gamma(x) \cdot \Gamma(y)}{\Gamma(x+y)}.
\end{equation*}
This result is also given in the supplement to \cite{E421}.

\subsection{The reflection formula}
Euler's version  of the reflection formula for the $\Gamma$-function, 
\begin{equation*}
\dfrac{\pi}{\sin \pi x}= \Gamma (x)\Gamma (1-x),
\end{equation*}
can be found in  $\S 43$ of \cite{E421} and reads
\[ [\lambda]\cdot [-\lambda]=\frac{\pi \lambda}{\sin \pi \lambda} ,\]
where $[\lambda]$ stands for $\lambda !$, that is $\Gamma(1+\lambda)$.\\
If one applies the reflection formula for $x=\frac{i}{n}$, $i=1,2,\cdots,n-1$, we obtain
$$
\begin{array}{rcl}
\Gamma \left(\dfrac{1}{n}\right)  \Gamma \left(\dfrac{n-1}{n}\right)&=&\dfrac{\pi}{\sin \frac{\pi}{n}},\\
\Gamma \left(\dfrac{2}{n}\right)  \Gamma \left(\dfrac{n-2}{n}\right)&=&\dfrac{\pi}{\sin \frac{2\pi}{n}},\\
\Gamma \left(\dfrac{3}{n}\right)  \Gamma \left(\dfrac{n-3}{n}\right)&=&\dfrac{\pi}{\sin \frac{3\pi}{n}},\\
\ldots&=&\ldots\\
\Gamma \left(\dfrac{n-1}{n}\right)  \Gamma \left(\dfrac{1}{n}\right)&=&\dfrac{\pi}{\sin \frac{(n-1)\pi}{n}} .\\
\end{array}
$$
Multiplying these equations together gives our first auxiliary formula

\begin{equation*}
 \prod_{i=1}^{n-1}\Gamma \left(\frac{i}{n}\right)^2= \frac{\pi^{n-1}}{\prod_{i=1}^{n-1} \sin \left(\frac{i \pi}{n}\right)} .
\end{equation*}

Our second auxiliary formula is
\begin{equation*}
\prod_{i=1}^{n-1} \sin \left(\frac{i \pi}{n} \right) = \frac{n}{2^{n-1}} ,
\end{equation*}
which is a nice exercise and which was certainly known to Euler.
For example, in $\S 7$ of \cite{E562} and in $\S 240$ of \cite{E101}, he states the more general formula

\begin{equation*}
\sin n \varphi = 2^{n-1} \sin \varphi \sin \left(\dfrac{\pi}{n}- \varphi\right) \sin \left(\dfrac{\pi}{n}+ \varphi\right)
\end{equation*}
\begin{equation*}
 \sin \left(\dfrac{2\pi}{n}- \varphi\right) \sin \left(\dfrac{2\pi}{n}+ \varphi\right)\cdot\text{etc.}
\end{equation*}
The product has $n$ factors in total. If we divide by $2^{n-1}\sin \varphi$, 
use $\sin \left(\frac{\pi(n-i)}{n}\right)= \sin \left(\frac{i\pi}{n}\right)$
and take the limit $\varphi \rightarrow 0$, we obtain the second auxiliary 
formula.\\[2mm]
The first and the second auxiliary formula were also given by Gauss in \cite{Ga28} and are used in his proof of the multiplication formula. Combining them and taking the square root, we obtain the beautiful formula
\begin{equation}
\Gamma \left(\frac{1}{n}\right) \Gamma \left( \frac{2}{n}\right) \cdots \Gamma \left(\frac{n-1}{n}\right)=\sqrt{\frac{(2\pi)^{n-1}}{n}}.
\end{equation}
This formula was also found by Euler in $\S 46$ of \cite{E816}, where he states it in the form

\begin{equation*}
\int\limits_{0}^{1}dx \left(\log \dfrac{1}{x}\right)^{\frac{1}{n}}\int\limits_{0}^{1}dx \left(\log \dfrac{1}{x}\right)^{\frac{2}{n}}\cdots \int\limits_{0}^{1}dx \left(\log \dfrac{1}{x}\right)^{\frac{n-1}{n}}= \dfrac{1 \cdot 2 \cdot 3 \cdots (n-1)}{n^{n-1}}\sqrt{\dfrac{2^{n-1}\pi^{n-1}}{n}}.
\end{equation*}

\section{Euler's version of the Multiplication Formula}

In $\S 53$ of \cite{E421} Euler gives the formula

\begin{equation*}
\left[ \frac{m}{n} \right] = \frac{m}{n} \sqrt[n]{n^{n-m}\cdot 1 \cdot 2 \cdot 3 \cdots (m-1) \left(\frac{1}{m}\right)\left(\frac{2}{m}\right)\left(\frac{3}{m}\right)\cdots \left(\frac{n-1}{m}\right)}.
\end{equation*}
As before, $[\lambda]$ is Euler's notation for the factorial of $\lambda$, so that $\left[ \frac{m}{n} \right] = \Gamma \left(\frac{m}{n}+1\right)$. Euler assumes $m$ and $n$ to be natural numbers, but it is easily seen that we can interpolate $1 \cdot 2 \cdot 3 \cdots (m-1)$ by $\Gamma(m)$. Therefore, if we assume $x$ to be real and positive and write $x$ instead of $m$ in the above formula and
express it in terms of the Beta-function using (2), Euler's formula becomes

\begin{equation*}
\Gamma \left(\frac{x}{n}\right)= \sqrt[n]{n^{n-x} \Gamma(x) \frac{1}{n^{n-1}}B \left(\frac{1}{n},\frac{x}{n}\right)B \left(\frac{2}{n},\frac{x}{n}\right)\cdots B \left(\frac{n-1}{n},\frac{x}{n}\right)}.
\end{equation*}
Expressing the Beta-function in terms of the $\Gamma$-function, then after 
some rearrangement under the $\sqrt[n]{}$-sign we obtain

\begin{equation*}
\Gamma \left(\frac{x}{n}\right)= \sqrt[n]{n^{1-x}\Gamma(x) \dfrac{\Gamma\left(\frac{1}{n}\right)\Gamma \left(\frac{x}{n}\right)}{\Gamma \left(\frac{x+1}{n}\right)} \cdot \dfrac{\Gamma\left(\frac{2}{n}\right)\Gamma \left(\frac{x}{n}\right)}{\Gamma \left(\frac{x+2}{n}\right)} \cdots \dfrac{\Gamma\left(\frac{n-1}{n}\right)\Gamma \left(\frac{x}{n}\right)}{\Gamma \left(\frac{x+n-1}{n}\right)}}.
\end{equation*}
By bringing all $\Gamma$-functions of fractional argument to the left-hand side,
the expression simplifies to

\begin{equation*}
\Gamma \left(\frac{x}{n}\right)\Gamma \left(\frac{x+1}{n}\right) \Gamma \left(\frac{x+2}{n}\right) \cdots \Gamma \left(\frac{x+n-1}{n}\right)= n^{1-x} \Gamma (x) \Gamma \left(\frac{1}{n}\right) \cdots \Gamma \left(\frac{n-1}{n}\right).
\end{equation*}
The product on the right-hand side, $\Gamma \left(\frac{1}{n}\right) \cdots \Gamma \left(\frac{n-1}{n}\right)$, was evaluated in (3) and thus we obtain 
\begin{equation*}
\Gamma \left(\frac{x}{n}\right)\Gamma \left(\frac{x+1}{n}\right) \Gamma \left(\frac{x+2}{n}\right) \cdots \Gamma \left(\frac{x+n-1}{n}\right)=  n^{1-x} \Gamma(x) \sqrt{\dfrac{(2\pi)^{n-1}}{n}}.
\end{equation*}
Thus, we arrived at the multiplication formula (1).

\section{Summary and Conclusion}

From the above sketch it is apparent that in \cite{E421} Euler had a result that is essentially equivalent to the 
multiplication formula for the $\Gamma$-function. He expressed it in terms of the symbol $\left(\frac{p}{q}\right)$, 
which is in modern notation is the Beta-function. One may wonder why Euler did not express his result in terms of the 
$\Gamma$-function itself. Reading his paper it becomes clear that his main motivation was to express the factorial of rational 
numbers in terms of integrals of \textit{algebraic functions}, and the formula 
given by Euler fulfills this purpose. For the same reason he probably did not replace $1 \cdot 2 \cdot 3 \cdots (m-1)$ by $\Gamma(m)$.\\

Euler also expressed $\Gamma(\frac{p}{q})$ in terms of integrals of algebraic algebraic functions in $\S 23$ \cite{E19} and $\S 5$ of \cite{E122}. That formula reads

\begin{equation*}
\int\limits_{0}^{1} \left(-\log x\right)^{\frac{p}{q}}dx =\sqrt[q]{1 \cdot 2 \cdot 3 \cdots p\left(\dfrac{2p}{q}+1\right)\left(\dfrac{3p}{q}+1\right)\left(\dfrac{4p}{q}+1\right)\cdots \left(\dfrac{qp}{q}+1\right)}
\end{equation*}
\begin{equation*}
\times \sqrt[q]{\int\limits_{0}^{1} dx(x-xx)^{\frac{p}{q}} \cdot \int\limits_{0}^{1} dx(x^2-x^3)^{\frac{p}{q}} \cdot \int\limits_{0}^{1} dx(x^3-x^4)^{\frac{p}{q}} \cdot \int\limits_{0}^{1} dx(x^4-x^5)^{\frac{p}{q}} \cdots \int\limits_{0}^{1} dx(x^{q-1}-x^q)^{\frac{p}{q}}}.
\end{equation*}
Despite the similarity to the first formula of section $2$, this formula is not as general as the multiplication formula\footnote{We want to mention here that in the foreword of the Opera Omnia, series 1, volume 19, p. LXI A. Krazer and G. Faber claim that these two formulas are equivalent and both are a special case of the multiplication formula. This is incorrect, as it was shown in the preceding sections. The formula given in section 3 does not lead to the multiplication formula, it only interpolates $\Gamma \left(\frac{p}{q}\right)$ in terms of algebraic integrals.}.
It appears that Euler was aware that the proofs he indicated in \cite{E421} were not completely convincing. 
He expressed that  with characteristic honesty in a concluding SCHOLIUM:\\

{\em  Hence infinitely many relations among the integral formulas of the form
\[ \int \dfrac{x^{p-1}dx}{(1-x^n)^{\frac{n-q}{n}}}=\left(\frac{p}{q} \right)\]
follow, which are even more remarkable, because we were led to them by a
completely singular method. And if anyone does not believe them to be true,
he or she should consult my observations on these integral formulas\footnote{Here Euler refers to his paper \cite{E321}.} 
and will then hence easily be convinced of their truth for any case. But even if
this consideration provides some confirmation, the relations found here are
nevertheless of even greater importance, because a certain structure is noticed
in them and they are easily generalized to all classes, whatever number was
assumed for the exponent $n$, whereas in the first treatment the calculation for the higher classes becomes continuously more cumbersome and intricate.}\\[2mm]
The history of the $\Gamma$-function is long and complex and apparently not all parts of the story have been 
told. We hope that this note provides some motivation to go carefully through other papers by Euler and other mathematicians of the past, not only because they make a good reading, but also to find further results, maybe stated in unfamiliar form, that were  proven rigorously by their successors. This will certainly be of interest for anyone studying the history 
of mathematics.

\section{Acknowledgements}

This author of this paper is supported by the \textit{Euler-Kreis Mainz.} The author especially wants to 
thank Prof. D. van Straten, Johannes Gutenberg Universität Mainz, both for very helpful suggestions 
concerning the presentation of the subject and for proof-reading the text.

\end{document}